\documentclass{article}

\usepackage{arxiv}

\usepackage[utf8]{inputenc}
\usepackage[T1]{fontenc}
\usepackage{hyperref}
\usepackage{url}
\usepackage{booktabs}
\usepackage{amsfonts}
\usepackage{nicefrac}
\usepackage{microtype}
\usepackage{lipsum}
\usepackage{algorithm, algpseudocode}
\usepackage{float}
\usepackage{tikz}
\usetikzlibrary{arrows}
\usepackage{capt-of}

\title{On the number of inequivalent \\ monotone Boolean functions of 8 variables}

\author{
  Bartłomiej Pawelski\\
    Institute of Informatics \\
    Faculty of Mathematics, Physics, and Informatics \\
  University of Gdańsk, 80-502 Gdańsk, Poland \\
  \texttt{bartlomiej.pawelski@ug.edu.pl} \\
}

\begin{document}
\maketitle

\vspace{-5mm}

\begin{abstract}
In this paper, the author presents algorithms that allow determining the number of fixed points in permutations of a set of monotone Boolean functions. Then, using Burnside's lemma, the author determines the number of inequivalent monotone Boolean functions of 8 variables. The number obtained is 1,392,195,548,889,993,358.
\end{abstract}

\vspace{2mm}

\keywords{Monotone Boolean functions \and Inequivalent monotone Boolean functions \and Dedekind numbers}

\section{Introduction}

Let $I^n$ be a set of all $n$ input variables of a Boolean function. Monotone Boolean function (MBF) is any Boolean function that can be implemented using only conjunctions and disjunctions \cite{r7}. Let $D_n$ be a set of all MBFs of $n$ variables, and $d_n$ a cardinality of this set; $d_n$ is also known as n-th Dedekind number.

Let $B^n$ be a power set of $I^n$. Each element in $B^n$ represents one of $2^n$ possible inputs of Boolean function. I am going to use notation in form: $\oslash$, $x_1$, $x_2$, $x_1 x_2$, $x_3$ ..., $x_1 x_2 x_3 ... x_n$ to describe elements in $B^n$.

I shall represent each Boolean function by a binary string of length $2^n$. Each $i$-th bit of function in this representation is Boolean output where the argument is an element from $B^n$ standing on the same position.

For example, consider the following truth table:

\begin{table}[H]
\centering
{\renewcommand{\arraystretch}{1.2}
\begin{tabular}{c|c|c|c|c|c|c|c}
$\oslash$ & $x_1$ & $x_2$ & $x_1 x_2$ & $x_3$ & $x_1 x_3$ & $x_2 x_3$ & $x_1 x_2 x_3$  \\ \hline
0 & 0 & 0 & 0 & 1 & 1 & 1 & 1 \\
\end{tabular}}
\caption{Example of MBF of three variables that returns true iff $x3$ is true} \label{tab:table1}
\end{table}

MBF from Table \ref{tab:table1} can be written as integer 15 for more convenient computer processing. All 6 MBFs in $D_2$ written as integers are: 0, 1, 3, 5, 7 and 15.

Two Boolean functions are equivalent if the first function can be transformed into the second function by any permutation of input variables. Let $R_n$ be a set of all inequivalent classes of $D_n$ (represented by lexicographically smallest representative MBF) and $r_n$ a cardinality of this set.

Any permutation of $I^n$ regroups elements of $B^n$ and $D_n$. There are $n!$ possible permutations of $I^n$, therefore there are at most $n!$ MBFs in one equivalence class.

In 1985, Chuchang and Shoben \cite{Ch1985} came up with the idea to derive a $r_n$ using the Burnside's Lemma. In the following year they calculated $r_7$ \cite{Ch1986}. Their result was confirmed by Stephen and Yusun in 2012 \cite{r7}. In 2018, Assarpour \cite{assarpour} has given lower bound of $r_8$ - 1,392,123,939,633,987,512.

In 1990, Wiedemann has calculated $d_8$ \cite{Wied1991}. He reduced the amount of work using properties of equivalent monotone Boolean functions. His result was confirmed in 2001 by Fidytek et al. \cite{Fid2001}.

In this paper I use Burnside's Lemma to derive $r_8$ = 1,392,195,548,889,993.358. I also describe the methods I developed for counting fixed points in $D_n$ under given permutation of $I^n$.

\begin{table}[H]
\centering
{\renewcommand{\arraystretch}{1.1}
\begin{tabular}{l|l|l}
$n$ & $d_n$                        & $r_n$                     \\ \hline
0 & 2                              & 2                         \\ 
1 & 3                              & 3                         \\ 
2 & 6                              & 5                         \\ 
3 & 20                             & 10                        \\ 
4 & 168                            & 30                        \\ 
5 & 7,581                          & 210                       \\ 
6 & 7,828,354                      & 16,353                    \\ 
7 & 2,414,682,040,998              & 490,013,148               \\
8 & 56,130,437,228,687,557,907,788 & 1,392,195,548,889,993,358 \\ \hline
\end{tabular}}
\caption{Known values of $d_n$ and $r_n$} \label{tab:values}
\end{table}

\section{Idea of counting IMBFs using Burnside's Lemma}

Chuchang and Shoben \cite{Ch1985} used the following application of Burnside's Lemma to count IMBFs of $n$ variables:

\begin{equation}
r_n = \frac{1}{n!} \sum_{i=1}^{k} \mu_i \phi(\pi_i)
\end{equation}

where:
\begin{itemize}
  \item $i$ = index of cycle type
  \item $\mu_i$ = number of permutations of certain cycle type $i$
  \item $\phi(\pi_i)$ = number of fixed points in $D_n$ under given permutation of $I^n$
\end{itemize}

The formula for determine $\mu_i$ for each cycle type is as follows:

\begin{equation}
\mu_i = \frac{n!}{(l^{k_1}_1 \times l^{k_2}_2 \cdot \cdot \cdot \times l^{k_r}_r)(k_1! \times k_2! \cdot \cdot \cdot \times k_r!)}
\end{equation}

with a cycle type of $k_1$ cycles of length $l_1$, $k_2$ cycles of length $l_2$, ... , $k_r$ cycles of length $l_r$  \cite{earl2014}. Precomputed results can be found in A181897 OEIS sequence.

The most difficult subproblem is fast counting fixed points in $D_n$ under a given permutation of $I^n$. The number of cycle types for the appropriate value of $n$ is described by the A000041 OEIS sequence. For $n=7$ there is 15 cycle types, and for $n=8$ there is 22 cycle types.

Let $\pi_i$ be the smallest representative cyclic permutation of $I^n$ with cycle type $i$.

\section{Algorithms counting fixed points in permutations}

For counting fixed points in $D_n$ after acting with specific permutation in the first place it is necessary to convert permutation of $n$ input variables into a permutation of $2^n$ bits in $B^n$ (hereinafter defined as $\pi_i B^n$). For example, one more time consider permutation of cycle type (123) and look how it regroups elements of $B^3$:

\begin{table}[H]
\centering
{\renewcommand{\arraystretch}{1.2}
\begin{tabular}{c|c|c|c|c|c|c|c|c}
& 0 & 1 & 2 & 3 & 4 & 5 & 6 & 7 \\
\hline
$(1)$ & $\oslash$ & $x_1$ & $x_2$ & $x_1 x_2$ & $x_3$ & $x_1 x_3$ & $x_2 x_3$ & $x_1 x_2 x_3$  \\
\hline
$(123)$ & $\oslash$ & $x_3$ & $x_1$ & $x_1 x_3$ & $x_2$ & $x_2 x_3$ & $x_1 x_2$ & $x_1 x_2 x_3$  \\
\end{tabular}}
\caption{Regroup of elements in $B^3$ under $\pi = (x_1\;x_2\;x_3)$} \label{tab:table4}
\end{table}

So $\pi = (x_1\;x_2\;x_3)$ should be transformed into $\pi(0)(124)(365)(7)B^3$. Each cycle designate points belonging to the same orbit. Points in each orbit must be set to the same value for the MBF to be fixed point of $D_3$ under $\pi = (x_1\;x_2\;x_3)$.

In this case, each element in $D_3$ to be a fixed point of $D_3$ under $\pi = (x_1\;x_2\;x_3)$ has to have:

\begin{itemize}
  \item 1-st, 2-nd and 4-th bit set on the same value
  \item 3-rd, 5-th and 6-th bit set on the same value
\end{itemize}

Using this approach fixed points in $D_3$ under $\pi = (x_1\;x_2\;x_3)$ can be simply found by iteration through all 20 elements in $D_3$ and checking which MBFs are satisfying the above conditions:

\begin{table}[H]
\centering
{\renewcommand{\arraystretch}{1.2}
\begin{tabular}{|c|c|c|c|c|c|c|c|c|}
\cline{2-9}
\multicolumn{1}{c}{} & \multicolumn{8}{|c|}{$n$-th bit of MBF} \\
\cline{1-9}
MBF written as integer & 0 & 1 & 2 & 3 & 4 & 5 & 6 & 7 \\
\hline
0 & 0 & 0 & 0 & 0 & 0 & 0 & 0 & 0  \\
\hline
1 & 0 & 0 & 0 & 0 & 0 & 0 & 0 & 1   \\
\hline
23 & 0 & 0 & 0 & 1 & 0 & 1 & 1 & 1   \\
\hline
127 & 0 & 1 & 1 & 1 & 1 & 1 & 1 & 1   \\
\hline
255 & 1 & 1 & 1 & 1 & 1 & 1 & 1 & 1    \\
\hline
\end{tabular}}
\caption{List of five fixed points in $D_3$ under $\pi = (x_1\;x_2\;x_3)$} \label{tab:table5}
\end{table}

\subsection{Generating a set of all fixed points in $D_n$ under permutation of cycle type of total length $n$}

\def\drawconnections{\draw (0,0) -- (-2,2);
    \draw (0,0) -- (2,2);
    \draw (-2,2) -- (-3,4);
    \draw (-2,2) -- (-1,4);
    \draw (-2,2) -- (1,4);
    \draw (2,2) -- (1,4);
    \draw (2,2) -- (3,4);
    \draw (2,2) -- (-1,4);
    \draw (-3,4) -- (-2,6);
    \draw (-1,4) -- (-2,6);
    \draw (-1,4) -- (2,6);
    \draw (1,4) -- (2,6);
    \draw (1,4) -- (-2,6);
    \draw (3,4) -- (2,6);
    \draw (-2,6) -- (0,8);
    \draw (2,6) -- (0,8);}

\def\dospace{\hspace{0.684cm}}

Instead of doing a lookup in $D_n$ for functions that have bits set the same value in each orbit, it is possible to generate a set of fixed points directly. I shall use fact that $D_n$ is equivalent to the set of all downsets of $B^n$ - so each element in $D_n$ is equivalent to some downset of $B^n$ \cite{campo}.

Two conditions must be met to generate MBF which is the fixed point in $D_n$ under the given permutation:

\begin{itemize}
  \item All points that are in one orbit should be set to the same value - 0 or 1.
  \item Value of points must respect the order of set inclusion.
\end{itemize}

For example, consider permutation $\pi = (x_1\;x_2)(x_3\;x_4)$. After transforming it into permutation of $B^4$, I get:

$\pi = $ (0)(1 2)(3)(4 8)(5 10)(6 9)(7 11)(12)(13 14)(15)$B^4$.

Now, let's transform this permutation into a poset of orbits, ordered by set inclusion. I shall represent orbits by their lexicographically smallest representative:

\begin{center}
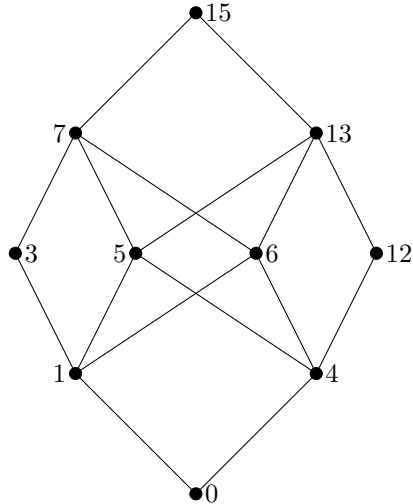


    \begin{tikzpicture}[scale=0.8]

    \drawconnections{}

    \draw[fill] (0,0) circle (.1cm) node[right] {$0$};
    \draw[fill] (-2,2) circle (.1cm) node[left] {$1$};
    \draw[fill] (-3,4) circle (.1cm) node[right] {$3$};
    \draw[fill] (2,2) circle (.1cm) node[right] {$4$};
    \draw[fill] (-1,4) circle (.1cm) node[left] {$5$};
    \draw[fill] (1,4) circle (.1cm) node[right] {$6$};
    \draw[fill] (-2,6) circle (.1cm) node[left] {$7$};
    \draw[fill] (3,4) circle (.1cm) node[right] {$12$};
    \draw[fill] (2,6) circle (.1cm) node[right] {$13$};
    \draw[fill] (0,8) circle (.1cm) node[right] {$15$};

    \end{tikzpicture}
    \captionof{figure}{Poset of orbits of $B_4$ under $\pi = (x_1\;x_2)(x_3\;x_4)$ ordered by set inclusion.}
    
\end{center}

Now it is only necessary to generate all downsets of this poset. In this case, the number of all downsets is 28:

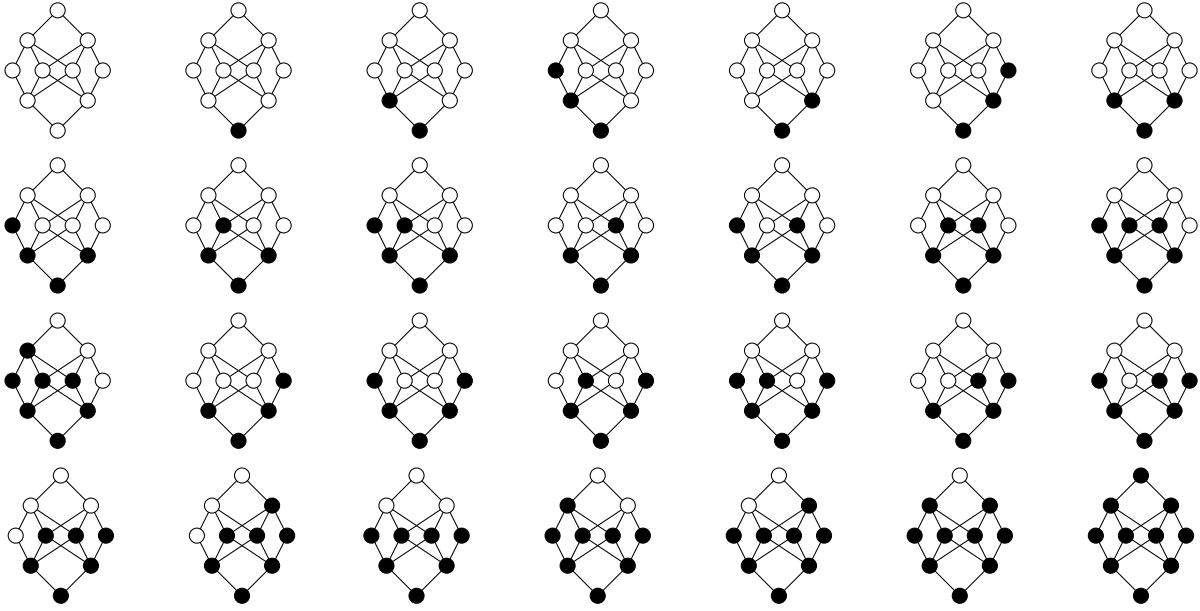
\vspace{0.7cm}
    \begin{center}
    \begin{tikzpicture}[scale=0.2]
    \drawconnections{}
    \draw[black, fill=white] (0,0) circle (.5cm) node[right] {};
    \draw[black, fill=white] (-2,2) circle (.5cm) node[right] {};
    \draw[black, fill=white] (-3,4) circle (.5cm) node[right] {};
    \draw[black, fill=white] (2,2) circle (.5cm) node[right] {};
    \draw[black, fill=white] (-1,4) circle (.5cm) node[left] {};
    \draw[black, fill=white] (1,4) circle (.5cm) node[right] {};
    \draw[black, fill=white] (-2,6) circle (.5cm) node[right] {};
    \draw[black, fill=white] (3,4) circle (.5cm) node[right] {};
    \draw[black, fill=white] (2,6) circle (.5cm) node[right] {};
    \draw[black, fill=white] (0,8) circle (.5cm) node[right] {};
    
    \end{tikzpicture} \dospace{} \begin{tikzpicture}[scale=0.2]
    \drawconnections{}
    \draw[black, fill=black] (0,0) circle (.5cm) node[right] {};
    \draw[black, fill=white] (-2,2) circle (.5cm) node[right] {};
    \draw[black, fill=white] (-3,4) circle (.5cm) node[right] {};
    \draw[black, fill=white] (2,2) circle (.5cm) node[right] {};
    \draw[black, fill=white] (-1,4) circle (.5cm) node[left] {};
    \draw[black, fill=white] (1,4) circle (.5cm) node[right] {};
    \draw[black, fill=white] (-2,6) circle (.5cm) node[right] {};
    \draw[black, fill=white] (3,4) circle (.5cm) node[right] {};
    \draw[black, fill=white] (2,6) circle (.5cm) node[right] {};
    \draw[black, fill=white] (0,8) circle (.5cm) node[right] {};

    \end{tikzpicture} \dospace{} \begin{tikzpicture}[scale=0.2]
    \drawconnections{}
    \draw[black, fill=black] (0,0) circle (.5cm) node[right] {};
    \draw[black, fill=black] (-2,2) circle (.5cm) node[right] {};
    \draw[black, fill=white] (-3,4) circle (.5cm) node[right] {};
    \draw[black, fill=white] (2,2) circle (.5cm) node[right] {};
    \draw[black, fill=white] (-1,4) circle (.5cm) node[left] {};
    \draw[black, fill=white] (1,4) circle (.5cm) node[right] {};
    \draw[black, fill=white] (-2,6) circle (.5cm) node[right] {};
    \draw[black, fill=white] (3,4) circle (.5cm) node[right] {};
    \draw[black, fill=white] (2,6) circle (.5cm) node[right] {};
    \draw[black, fill=white] (0,8) circle (.5cm) node[right] {};

    \end{tikzpicture} \dospace{} \begin{tikzpicture}[scale=0.2]
    \drawconnections{}
    \draw[black, fill=black] (0,0) circle (.5cm) node[right] {};
    \draw[black, fill=black] (-2,2) circle (.5cm) node[right] {};
    \draw[black, fill=black] (-3,4) circle (.5cm) node[right] {};
    \draw[black, fill=white] (2,2) circle (.5cm) node[right] {};
    \draw[black, fill=white] (-1,4) circle (.5cm) node[left] {};
    \draw[black, fill=white] (1,4) circle (.5cm) node[right] {};
    \draw[black, fill=white] (-2,6) circle (.5cm) node[right] {};
    \draw[black, fill=white] (3,4) circle (.5cm) node[right] {};
    \draw[black, fill=white] (2,6) circle (.5cm) node[right] {};
    \draw[black, fill=white] (0,8) circle (.5cm) node[right] {};
    \end{tikzpicture} \dospace{} \begin{tikzpicture}[scale=0.2]
    \drawconnections{}
    \draw[black, fill=black] (0,0) circle (.5cm) node[right] {};
    \draw[black, fill=white] (-2,2) circle (.5cm) node[right] {};
    \draw[black, fill=white] (-3,4) circle (.5cm) node[right] {};
    \draw[black, fill=black] (2,2) circle (.5cm) node[right] {};
    \draw[black, fill=white] (-1,4) circle (.5cm) node[left] {};
    \draw[black, fill=white] (1,4) circle (.5cm) node[right] {};
    \draw[black, fill=white] (-2,6) circle (.5cm) node[right] {};
    \draw[black, fill=white] (3,4) circle (.5cm) node[right] {};
    \draw[black, fill=white] (2,6) circle (.5cm) node[right] {};
    \draw[black, fill=white] (0,8) circle (.5cm) node[right] {};
    \end{tikzpicture} \dospace{} \begin{tikzpicture}[scale=0.2]
    \drawconnections{}
    \draw[black, fill=black] (0,0) circle (.5cm) node[right] {};
    \draw[black, fill=white] (-2,2) circle (.5cm) node[right] {};
    \draw[black, fill=white] (-3,4) circle (.5cm) node[right] {};
    \draw[black, fill=black] (2,2) circle (.5cm) node[right] {};
    \draw[black, fill=white] (-1,4) circle (.5cm) node[left] {};
    \draw[black, fill=white] (1,4) circle (.5cm) node[right] {};
    \draw[black, fill=white] (-2,6) circle (.5cm) node[right] {};
    \draw[black, fill=black] (3,4) circle (.5cm) node[right] {};
    \draw[black, fill=white] (2,6) circle (.5cm) node[right] {};
    \draw[black, fill=white] (0,8) circle (.5cm) node[right] {};
    \end{tikzpicture} \dospace{} \begin{tikzpicture}[scale=0.2]
    \drawconnections{}
    \draw[black, fill=black] (0,0) circle (.5cm) node[right] {};
    \draw[black, fill=black] (-2,2) circle (.5cm) node[right] {};
    \draw[black, fill=white] (-3,4) circle (.5cm) node[right] {};
    \draw[black, fill=black] (2,2) circle (.5cm) node[right] {};
    \draw[black, fill=white] (-1,4) circle (.5cm) node[left] {};
    \draw[black, fill=white] (1,4) circle (.5cm) node[right] {};
    \draw[black, fill=white] (-2,6) circle (.5cm) node[right] {};
    \draw[black, fill=white] (3,4) circle (.5cm) node[right] {};
    \draw[black, fill=white] (2,6) circle (.5cm) node[right] {};
    \draw[black, fill=white] (0,8) circle (.5cm) node[right] {};
    \end{tikzpicture}
    
    \begin{tikzpicture}[scale=0.2]
    \drawconnections{}
    \draw[black, fill=black] (0,0) circle (.5cm) node[right] {};
    \draw[black, fill=black] (-2,2) circle (.5cm) node[right] {};
    \draw[black, fill=black] (-3,4) circle (.5cm) node[right] {};
    \draw[black, fill=black] (2,2) circle (.5cm) node[right] {};
    \draw[black, fill=white] (-1,4) circle (.5cm) node[left] {};
    \draw[black, fill=white] (1,4) circle (.5cm) node[right] {};
    \draw[black, fill=white] (-2,6) circle (.5cm) node[right] {};
    \draw[black, fill=white] (3,4) circle (.5cm) node[right] {};
    \draw[black, fill=white] (2,6) circle (.5cm) node[right] {};
    \draw[black, fill=white] (0,8) circle (.5cm) node[right] {};
    \end{tikzpicture} \dospace{} \begin{tikzpicture}[scale=0.2]
    \drawconnections{}
    \draw[black, fill=black] (0,0) circle (.5cm) node[right] {};
    \draw[black, fill=black] (-2,2) circle (.5cm) node[right] {};
    \draw[black, fill=white] (-3,4) circle (.5cm) node[right] {};
    \draw[black, fill=black] (2,2) circle (.5cm) node[right] {};
    \draw[black, fill=black] (-1,4) circle (.5cm) node[left] {};
    \draw[black, fill=white] (1,4) circle (.5cm) node[right] {};
    \draw[black, fill=white] (-2,6) circle (.5cm) node[right] {};
    \draw[black, fill=white] (3,4) circle (.5cm) node[right] {};
    \draw[black, fill=white] (2,6) circle (.5cm) node[right] {};
    \draw[black, fill=white] (0,8) circle (.5cm) node[right] {};
    \end{tikzpicture} \dospace{} \begin{tikzpicture}[scale=0.2]
    \drawconnections{}
    \draw[black, fill=black] (0,0) circle (.5cm) node[right] {};
    \draw[black, fill=black] (-2,2) circle (.5cm) node[right] {};
    \draw[black, fill=black] (-3,4) circle (.5cm) node[right] {};
    \draw[black, fill=black] (2,2) circle (.5cm) node[right] {};
    \draw[black, fill=black] (-1,4) circle (.5cm) node[left] {};
    \draw[black, fill=white] (1,4) circle (.5cm) node[right] {};
    \draw[black, fill=white] (-2,6) circle (.5cm) node[right] {};
    \draw[black, fill=white] (3,4) circle (.5cm) node[right] {};
    \draw[black, fill=white] (2,6) circle (.5cm) node[right] {};
    \draw[black, fill=white] (0,8) circle (.5cm) node[right] {};
    \end{tikzpicture} \dospace{} \begin{tikzpicture}[scale=0.2]
    \drawconnections{}
    \draw[black, fill=black] (0,0) circle (.5cm) node[right] {};
    \draw[black, fill=black] (-2,2) circle (.5cm) node[right] {};
    \draw[black, fill=white] (-3,4) circle (.5cm) node[right] {};
    \draw[black, fill=black] (2,2) circle (.5cm) node[right] {};
    \draw[black, fill=white] (-1,4) circle (.5cm) node[left] {};
    \draw[black, fill=black] (1,4) circle (.5cm) node[right] {};
    \draw[black, fill=white] (-2,6) circle (.5cm) node[right] {};
    \draw[black, fill=white] (3,4) circle (.5cm) node[right] {};
    \draw[black, fill=white] (2,6) circle (.5cm) node[right] {};
    \draw[black, fill=white] (0,8) circle (.5cm) node[right] {};
    \end{tikzpicture} \dospace{} \begin{tikzpicture}[scale=0.2]
    \drawconnections{}
    \draw[black, fill=black] (0,0) circle (.5cm) node[right] {};
    \draw[black, fill=black] (-2,2) circle (.5cm) node[right] {};
    \draw[black, fill=black] (-3,4) circle (.5cm) node[right] {};
    \draw[black, fill=black] (2,2) circle (.5cm) node[right] {};
    \draw[black, fill=white] (-1,4) circle (.5cm) node[left] {};
    \draw[black, fill=black] (1,4) circle (.5cm) node[right] {};
    \draw[black, fill=white] (-2,6) circle (.5cm) node[right] {};
    \draw[black, fill=white] (3,4) circle (.5cm) node[right] {};
    \draw[black, fill=white] (2,6) circle (.5cm) node[right] {};
    \draw[black, fill=white] (0,8) circle (.5cm) node[right] {};
    \end{tikzpicture} \dospace{} \begin{tikzpicture}[scale=0.2]
    \drawconnections{}
    \draw[black, fill=black] (0,0) circle (.5cm) node[right] {};
    \draw[black, fill=black] (-2,2) circle (.5cm) node[right] {};
    \draw[black, fill=white] (-3,4) circle (.5cm) node[right] {};
    \draw[black, fill=black] (2,2) circle (.5cm) node[right] {};
    \draw[black, fill=black] (-1,4) circle (.5cm) node[left] {};
    \draw[black, fill=black] (1,4) circle (.5cm) node[right] {};
    \draw[black, fill=white] (-2,6) circle (.5cm) node[right] {};
    \draw[black, fill=white] (3,4) circle (.5cm) node[right] {};
    \draw[black, fill=white] (2,6) circle (.5cm) node[right] {};
    \draw[black, fill=white] (0,8) circle (.5cm) node[right] {};
    \end{tikzpicture} \dospace{} \begin{tikzpicture}[scale=0.2]
    \drawconnections{}
    \draw[black, fill=black] (0,0) circle (.5cm) node[right] {};
    \draw[black, fill=black] (-2,2) circle (.5cm) node[right] {};
    \draw[black, fill=black] (-3,4) circle (.5cm) node[right] {};
    \draw[black, fill=black] (2,2) circle (.5cm) node[right] {};
    \draw[black, fill=black] (-1,4) circle (.5cm) node[left] {};
    \draw[black, fill=black] (1,4) circle (.5cm) node[right] {};
    \draw[black, fill=white] (-2,6) circle (.5cm) node[right] {};
    \draw[black, fill=white] (3,4) circle (.5cm) node[right] {};
    \draw[black, fill=white] (2,6) circle (.5cm) node[right] {};
    \draw[black, fill=white] (0,8) circle (.5cm) node[right] {};
    \end{tikzpicture}
    
    \begin{tikzpicture}[scale=0.2]
    \drawconnections{}
    \draw[black, fill=black] (0,0) circle (.5cm) node[right] {};
    \draw[black, fill=black] (-2,2) circle (.5cm) node[right] {};
    \draw[black, fill=black] (-3,4) circle (.5cm) node[right] {};
    \draw[black, fill=black] (2,2) circle (.5cm) node[right] {};
    \draw[black, fill=black] (-1,4) circle (.5cm) node[left] {};
    \draw[black, fill=black] (1,4) circle (.5cm) node[right] {};
    \draw[black, fill=black] (-2,6) circle (.5cm) node[right] {};
    \draw[black, fill=white] (3,4) circle (.5cm) node[right] {};
    \draw[black, fill=white] (2,6) circle (.5cm) node[right] {};
    \draw[black, fill=white] (0,8) circle (.5cm) node[right] {};
    \end{tikzpicture} \dospace{} \begin{tikzpicture}[scale=0.2]
    \drawconnections{}
    \draw[black, fill=black] (0,0) circle (.5cm) node[right] {};
    \draw[black, fill=black] (-2,2) circle (.5cm) node[right] {};
    \draw[black, fill=white] (-3,4) circle (.5cm) node[right] {};
    \draw[black, fill=black] (2,2) circle (.5cm) node[right] {};
    \draw[black, fill=white] (-1,4) circle (.5cm) node[left] {};
    \draw[black, fill=white] (1,4) circle (.5cm) node[right] {};
    \draw[black, fill=white] (-2,6) circle (.5cm) node[right] {};
    \draw[black, fill=black] (3,4) circle (.5cm) node[right] {};
    \draw[black, fill=white] (2,6) circle (.5cm) node[right] {};
    \draw[black, fill=white] (0,8) circle (.5cm) node[right] {};
    \end{tikzpicture} \dospace{} \begin{tikzpicture}[scale=0.2]
    \drawconnections{}
    \draw[black, fill=black] (0,0) circle (.5cm) node[right] {};
    \draw[black, fill=black] (-2,2) circle (.5cm) node[right] {};
    \draw[black, fill=black] (-3,4) circle (.5cm) node[right] {};
    \draw[black, fill=black] (2,2) circle (.5cm) node[right] {};
    \draw[black, fill=white] (-1,4) circle (.5cm) node[left] {};
    \draw[black, fill=white] (1,4) circle (.5cm) node[right] {};
    \draw[black, fill=white] (-2,6) circle (.5cm) node[right] {};
    \draw[black, fill=black] (3,4) circle (.5cm) node[right] {};
    \draw[black, fill=white] (2,6) circle (.5cm) node[right] {};
    \draw[black, fill=white] (0,8) circle (.5cm) node[right] {};
    \end{tikzpicture} \dospace{} \begin{tikzpicture}[scale=0.2]
    \drawconnections{}
    \draw[black, fill=black] (0,0) circle (.5cm) node[right] {};
    \draw[black, fill=black] (-2,2) circle (.5cm) node[right] {};
    \draw[black, fill=white] (-3,4) circle (.5cm) node[right] {};
    \draw[black, fill=black] (2,2) circle (.5cm) node[right] {};
    \draw[black, fill=black] (-1,4) circle (.5cm) node[left] {};
    \draw[black, fill=white] (1,4) circle (.5cm) node[right] {};
    \draw[black, fill=white] (-2,6) circle (.5cm) node[right] {};
    \draw[black, fill=black] (3,4) circle (.5cm) node[right] {};
    \draw[black, fill=white] (2,6) circle (.5cm) node[right] {};
    \draw[black, fill=white] (0,8) circle (.5cm) node[right] {};
    \end{tikzpicture} \dospace{} \begin{tikzpicture}[scale=0.2]
    \drawconnections{}
    \draw[black, fill=black] (0,0) circle (.5cm) node[right] {};
    \draw[black, fill=black] (-2,2) circle (.5cm) node[right] {};
    \draw[black, fill=black] (-3,4) circle (.5cm) node[right] {};
    \draw[black, fill=black] (2,2) circle (.5cm) node[right] {};
    \draw[black, fill=black] (-1,4) circle (.5cm) node[left] {};
    \draw[black, fill=white] (1,4) circle (.5cm) node[right] {};
    \draw[black, fill=white] (-2,6) circle (.5cm) node[right] {};
    \draw[black, fill=black] (3,4) circle (.5cm) node[right] {};
    \draw[black, fill=white] (2,6) circle (.5cm) node[right] {};
    \draw[black, fill=white] (0,8) circle (.5cm) node[right] {};
    \end{tikzpicture} \dospace{} \begin{tikzpicture}[scale=0.2]
    \drawconnections{}
    \draw[black, fill=black] (0,0) circle (.5cm) node[right] {};
    \draw[black, fill=black] (-2,2) circle (.5cm) node[right] {};
    \draw[black, fill=white] (-3,4) circle (.5cm) node[right] {};
    \draw[black, fill=black] (2,2) circle (.5cm) node[right] {};
    \draw[black, fill=white] (-1,4) circle (.5cm) node[left] {};
    \draw[black, fill=black] (1,4) circle (.5cm) node[right] {};
    \draw[black, fill=white] (-2,6) circle (.5cm) node[right] {};
    \draw[black, fill=black] (3,4) circle (.5cm) node[right] {};
    \draw[black, fill=white] (2,6) circle (.5cm) node[right] {};
    \draw[black, fill=white] (0,8) circle (.5cm) node[right] {};
    \end{tikzpicture} \dospace{} \begin{tikzpicture}[scale=0.2]
    \drawconnections{}
    \draw[black, fill=black] (0,0) circle (.5cm) node[right] {};
    \draw[black, fill=black] (-2,2) circle (.5cm) node[right] {};
    \draw[black, fill=black] (-3,4) circle (.5cm) node[right] {};
    \draw[black, fill=black] (2,2) circle (.5cm) node[right] {};
    \draw[black, fill=white] (-1,4) circle (.5cm) node[left] {};
    \draw[black, fill=black] (1,4) circle (.5cm) node[right] {};
    \draw[black, fill=white] (-2,6) circle (.5cm) node[right] {};
    \draw[black, fill=black] (3,4) circle (.5cm) node[right] {};
    \draw[black, fill=white] (2,6) circle (.5cm) node[right] {};
    \draw[black, fill=white] (0,8) circle (.5cm) node[right] {};
    \end{tikzpicture}
    
    \begin{tikzpicture}[scale=0.2]
    \drawconnections{}
    \draw[black, fill=black] (0,0) circle (.5cm) node[right] {};
    \draw[black, fill=black] (-2,2) circle (.5cm) node[right] {};
    \draw[black, fill=white] (-3,4) circle (.5cm) node[right] {};
    \draw[black, fill=black] (2,2) circle (.5cm) node[right] {};
    \draw[black, fill=black] (-1,4) circle (.5cm) node[left] {};
    \draw[black, fill=black] (1,4) circle (.5cm) node[right] {};
    \draw[black, fill=white] (-2,6) circle (.5cm) node[right] {};
    \draw[black, fill=black] (3,4) circle (.5cm) node[right] {};
    \draw[black, fill=white] (2,6) circle (.5cm) node[right] {};
    \draw[black, fill=white] (0,8) circle (.5cm) node[right] {};
    \end{tikzpicture} \dospace{} \begin{tikzpicture}[scale=0.2]
    \drawconnections{}
    \draw[black, fill=black] (0,0) circle (.5cm) node[right] {};
    \draw[black, fill=black] (-2,2) circle (.5cm) node[right] {};
    \draw[black, fill=white] (-3,4) circle (.5cm) node[right] {};
    \draw[black, fill=black] (2,2) circle (.5cm) node[right] {};
    \draw[black, fill=black] (-1,4) circle (.5cm) node[left] {};
    \draw[black, fill=black] (1,4) circle (.5cm) node[right] {};
    \draw[black, fill=white] (-2,6) circle (.5cm) node[right] {};
    \draw[black, fill=black] (3,4) circle (.5cm) node[right] {};
    \draw[black, fill=black] (2,6) circle (.5cm) node[right] {};
    \draw[black, fill=white] (0,8) circle (.5cm) node[right] {};
    \end{tikzpicture} \dospace{}\begin{tikzpicture}[scale=0.2]
    \drawconnections{}
    \draw[black, fill=black] (0,0) circle (.5cm) node[right] {};
    \draw[black, fill=black] (-2,2) circle (.5cm) node[right] {};
    \draw[black, fill=black] (-3,4) circle (.5cm) node[right] {};
    \draw[black, fill=black] (2,2) circle (.5cm) node[right] {};
    \draw[black, fill=black] (-1,4) circle (.5cm) node[left] {};
    \draw[black, fill=black] (1,4) circle (.5cm) node[right] {};
    \draw[black, fill=white] (-2,6) circle (.5cm) node[right] {};
    \draw[black, fill=black] (3,4) circle (.5cm) node[right] {};
    \draw[black, fill=white] (2,6) circle (.5cm) node[right] {};
    \draw[black, fill=white] (0,8) circle (.5cm) node[right] {};
    \end{tikzpicture} \dospace{} \begin{tikzpicture}[scale=0.2]
    \drawconnections{}
    \draw[black, fill=black] (0,0) circle (.5cm) node[right] {};
    \draw[black, fill=black] (-2,2) circle (.5cm) node[right] {};
    \draw[black, fill=black] (-3,4) circle (.5cm) node[right] {};
    \draw[black, fill=black] (2,2) circle (.5cm) node[right] {};
    \draw[black, fill=black] (-1,4) circle (.5cm) node[left] {};
    \draw[black, fill=black] (1,4) circle (.5cm) node[right] {};
    \draw[black, fill=black] (-2,6) circle (.5cm) node[right] {};
    \draw[black, fill=black] (3,4) circle (.5cm) node[right] {};
    \draw[black, fill=white] (2,6) circle (.5cm) node[right] {};
    \draw[black, fill=white] (0,8) circle (.5cm) node[right] {};
    \end{tikzpicture} \dospace{} \begin{tikzpicture}[scale=0.2]
    \drawconnections{}
    \draw[black, fill=black] (0,0) circle (.5cm) node[right] {};
    \draw[black, fill=black] (-2,2) circle (.5cm) node[right] {};
    \draw[black, fill=black] (-3,4) circle (.5cm) node[right] {};
    \draw[black, fill=black] (2,2) circle (.5cm) node[right] {};
    \draw[black, fill=black] (-1,4) circle (.5cm) node[left] {};
    \draw[black, fill=black] (1,4) circle (.5cm) node[right] {};
    \draw[black, fill=white] (-2,6) circle (.5cm) node[right] {};
    \draw[black, fill=black] (3,4) circle (.5cm) node[right] {};
    \draw[black, fill=black] (2,6) circle (.5cm) node[right] {};
    \draw[black, fill=white] (0,8) circle (.5cm) node[right] {};
    \end{tikzpicture} \dospace{} \begin{tikzpicture}[scale=0.2]
    \drawconnections{}
    \draw[black, fill=black] (0,0) circle (.5cm) node[right] {};
    \draw[black, fill=black] (-2,2) circle (.5cm) node[right] {};
    \draw[black, fill=black] (-3,4) circle (.5cm) node[right] {};
    \draw[black, fill=black] (2,2) circle (.5cm) node[right] {};
    \draw[black, fill=black] (-1,4) circle (.5cm) node[left] {};
    \draw[black, fill=black] (1,4) circle (.5cm) node[right] {};
    \draw[black, fill=black] (-2,6) circle (.5cm) node[right] {};
    \draw[black, fill=black] (3,4) circle (.5cm) node[right] {};
    \draw[black, fill=black] (2,6) circle (.5cm) node[right] {};
    \draw[black, fill=white] (0,8) circle (.5cm) node[right] {};
    \end{tikzpicture} \dospace{} \begin{tikzpicture}[scale=0.2]
    \drawconnections{}
    \draw[black, fill=black] (0,0) circle (.5cm) node[right] {};
    \draw[black, fill=black] (-2,2) circle (.5cm) node[right] {};
    \draw[black, fill=black] (-3,4) circle (.5cm) node[right] {};
    \draw[black, fill=black] (2,2) circle (.5cm) node[right] {};
    \draw[black, fill=black] (-1,4) circle (.5cm) node[left] {};
    \draw[black, fill=black] (1,4) circle (.5cm) node[right] {};
    \draw[black, fill=black] (-2,6) circle (.5cm) node[right] {};
    \draw[black, fill=black] (3,4) circle (.5cm) node[right] {};
    \draw[black, fill=black] (2,6) circle (.5cm) node[right] {};
    \draw[black, fill=black] (0,8) circle (.5cm) node[right] {};
    \end{tikzpicture}

\end{center}
\vspace{1cm}
The structures thus obtained are equivalent to the set of fixed points in $D_4$ under $\pi = (x_1\;x_2)(x_3\;x_4)$. One can unpack the downsets obtained in this way to the form of MBF of $2^n$ length, but there is no need for that, as long as it is only interesting how many fixed points there are.

I am using this algorithm only to generate a set of fixes of permutation with cycle type of total length $n$ - for example, I use Algorithm 1 to create set of fixed points in $D_4$ under $\pi = (x_1\;x_2)(x_3\;x_4)$, but to create a set of $D_5$ under $\pi = (x_1\;x_2)(x_3\;x_4)$ (or higher) it is more convenient and cheaper computationally to use Algorithm 2.

\begin{algorithm}[H]
\caption{Generate a set of fixed points in $D_n$ under permutation of cycle type of total length $n$}
    \vspace{1mm}
    \hspace*{\algorithmicindent} \textbf{Input:} Cycle type $i$ of total length $n$ \\
    \hspace*{\algorithmicindent} \textbf{Output:} Set $S$ of all fixed points in $D_n$
\begin{algorithmic}[1]
\State Transform $\pi_iI^n$ into $\pi_iB^n$
\State Generate set $Orb_i$ containing all orbits in $\pi_iB^n$
\State Order $Orb_i$ into poset $P$ by set inclusion
\State Initialize set $S$ of downsets of $P$
\State Add two downsets: \{\} and \{0\} to $S$
\ForAll {\textbf{elements} $a \in P $}
\ForAll {\textbf{elements} $b \in S $}
\If{({b} $\cup$ a) is downset of $P$}
  \State Add downset ({b} $\cup$ a) to $S$
\EndIf
\EndFor
\EndFor
\end{algorithmic}
\end{algorithm}

\subsection{Generating a set of all fixed points in $D_{n+1}$ under permutation of cycle type of total length $n$}

Each MBF in $D_{n+1}$ consists of two concatenated functions ($\alpha$, $\beta$) from the set $D_n$. Moreover, there must be a relation $\alpha \preceq \beta$, which means that for every i-th bit $\alpha_i \leq \beta_i$ \cite{Bakoev}.

If MBF under permutation of cycle type $\pi_i$ of total length $n$ has $n$+$1$ variables, then $n$+$1$-th variable is fixed point of $\pi_i I^{n+1}$. This means that the first half ($\alpha$) of $B_{n+1}$ will be affected by the permutation $\pi_i$. The second half ($\beta$) of $B_{n+1}$ has to be in relation "$\preceq$" with the first half, which means that it will be have set at least the same bits. Moreover, the first $n$ variables are regrouped in the same way as in the first half and will therefore satisfy conditions of being a fixed point of $\pi_i D_n$.

The conclusion from this reasoning is as follows: each MBF from $D_{n+1}$ being a fixed point in $\pi_i$ of total length $n$ consist of two MBFs ($\alpha, \beta, \alpha \preceq \beta$) being a fixed points in $\pi_i D_n$.

\begin{table}[H]
\centering
\begin{tabular}{c|c|c|c|c|c|c|c|c}
& 0 & 1 & 2 & 3 & 4 & 5 & 6 & 7 \\
\hline
$id$ & $\oslash$ & $x_1$ & $x_2$ & $x_1 x_2$ & $x_3$ & $x_1 x_3$ & $x_2 x_3$ & $x_1 x_2 x_3$  \\
\hline
$(12)$ & $\oslash$ & \textbf{$x_2$} & $x_1$ & $x_1 x_2$ & $x_3$ & $x_2 x_3$ & $x_1 x_3$ & $x_1 x_2 x_3$  \\
\end{tabular}
\caption{Regroup of elements in $B^3$ under $\pi = (x_1\;x_2)$} \label{tab:regroupB3}
\end{table}

For example look on the Table \ref{tab:regroupB3}, where it is visible that first half of $B_3$ under $\pi = (x_1\;x_2)$ is the same as $B_2$ under $\pi = (x_1\;x_2)$. It also visible that $x_3$ is fixed point in $B_3$ under $\pi = (x_1\;x_2)$ - it is at the same places as in $B_3$.

Having a set of MBFs being a set of fixed points of $\pi_i D_n$ it is only necessary to use well-known algorithms used for determining Dedekind numbers (for example \cite{Wied1991}, \cite{Fid2001}), but instead of giving $D_n$ on input, set of fixed points in $\pi_i D_n$ will be given.

For example, let's determine number of fixed points in $D_5$ under $\pi = (x_1\;x_2)(x_3\;x_4)$. From previous subsection it is known that there are a 28 elements in the set of fixed points of $D_4$ under $\pi = (x_1\;x_2)(x_3\;x_4)$.

I will use modified "Algorithm 1" from \cite{Fid2001} (any algorithm from this paper will do the job, however, algorithms 2 and 3 don't give a set but its cardinality). For each pair of MBFs of length $2^n$ bits in this set, it must be checked if there is a relation "$\preceq$" between them. It is convenient to use bitwise OR for this task.

\begin{algorithm}[H]
\caption{Generating a set of all fixed points in $D_{n+1}$ under permutation of cycle type of total length $n$}
    \vspace{1mm}
    \hspace*{\algorithmicindent} \textbf{Input:} Cycle type $i$ of total length $n$ \\
    \hspace*{\algorithmicindent} \textbf{Output:} Set $S$ of all fixed points in $D_{n+1}$ written as integer
\begin{algorithmic}[1]
\State Use Algorithm 1 to generate a set $S'$ of all fixed points in $\pi_i D_n$
\State Convert all MBFs in $S'$ to integers of length $2^n$ bits
\State Initialize set $S$ of integers of length $2^{n+1}$ bits
\ForAll {\textbf{elements} $a \in S' $}
\ForAll {\textbf{elements} $b \in S' $}
\If{(a | b) = b} \Comment{"|" is bitwise "OR"}
  \State Add integer ((a << $2^n$) | b) to $S$ \Comment{"<<" is logical shift}
\EndIf
\EndFor
\EndFor
\end{algorithmic}
\end{algorithm}

\subsection{Special case: determining cardinality of set of fixed points of $D_8$ under $\pi = (x_1\;x_2)(x_3\;x_4)(x_5\;x_6)(x_7\;x_8)$}

Determining number of fixed points in $D_8$ under $\pi = (x_1\;x_2)(x_3\;x_4)(x_5\;x_6)(x_7\;x_8)$ was too memory-intensive for Algorithm 1 considering the resources at hand. The size of the largest antichain (or simply: width) of poset of orbits of super set of $\pi = (x_1\;x_2)(x_3\;x_4)(x_5\;x_6)(x_7\;x_8)$ is 38, so the weak lower bound of number of fixed points in $D_8$ under $\pi = (x_1\;x_2)(x_3\;x_4)(x_5\;x_6)(x_7\;x_8)$ = 274877906944. In practice, even the machine with 128GB RAM was insufficient - so there was a need to develop a better algorithm for this particular case.

The idea of a cheaper calculation of this number was based on noticing the following:

Look at Wiedemann's approach \cite{Wied1991}. He noticed that each MBF from the set $D_{n+2}$ can be splitted into 4 MBFs from $D_n$: $\alpha_w, \beta_w, \gamma_w, \delta_w$, and there are following dependencies:

\begin{itemize}
  \item $\alpha_w \preceq \beta_w$
  \item $\alpha_w \preceq \gamma_w$
  \item $\beta_w \preceq \delta_w$
  \item $\gamma_w \preceq \delta_w$
\end{itemize}

Each fixed point in $D_{n+2}$ under permutation of cycle type containing only 1-cycles and 2-cycles (let $\pi_{(i)(n+1\:n+2)}$ be definition of this permutation) can be splitted into 4 MBFs from $D_n$:

\begin{itemize}
  \item $\alpha$ = being MBF from the set of fixed points in $\pi_i D_n$
  \item $\beta$ = being MBF from the set $D_n$
  \item $\gamma$ = being MBF from the set $D_n$
  \item $\delta$ = being MBF from the set of fixed points in $\pi_i D_n$
\end{itemize}

And there are the following dependencies:

\begin{itemize}
  \item $\alpha \preceq \beta$
  \item $\alpha \preceq \gamma$
  \item $\beta \preceq \delta$
  \item $\gamma \preceq \delta$
\end{itemize}

And there are only $d_n$ valid pairs of $\beta \gamma$, because for each bit of these MBFs $\beta_i$ is equal to $\pi(12)\gamma_i$.

For example, cycle type of $B^4$ under $\pi = (x_1\;x_2)(x_3\;x_4)$ was $\pi$(0)(1 2)(3)(4 8)(5 10)(6 9)(7 11)(12)(13 14)(15)$B^4$. Let's break it down into three parts following the above approach:

\begin{itemize}
  \item $\alpha$ as (0)(1 2)(3); being MBF from the set of fixed points in $\pi = (x_1\;x_2) D_n$
  \item $\beta \gamma$ as (4 8)(5 10)(6 9)(7 11) being pairs of MBFs from the set $D_n$ such that for each bit $\beta_i$ = ($\pi = (x_1\;x_2))\gamma_i$
  \item $\delta$ as (12)(13 14)(15), being MBF from the set of fixed points in $\pi = (x_1\;x_2) D_n$
\end{itemize}

\begin{algorithm}[H]
\caption{Determining a cardinality of set of all fixed points in $D_{n+2}$ under $\pi_{(i)(n+1\:n+2)}$}
    \vspace{1mm}
    \hspace*{\algorithmicindent} \textbf{Input:} $D_n$ and set of fixed points in $\pi_i D_n$ with MBFs written as integers  \\
    \hspace*{\algorithmicindent} \textbf{Output:} $k$: cardinality of set of fixed points of $D_{n+2}$ under $\pi_{(i)(n+1\:n+2)}$
\begin{algorithmic}[1]
\State Initialize $k = 0$,
\ForAll {$\beta \in D_n $}
\State Determine $\gamma$ such that for each bit in ($\beta, \gamma$) \ $\beta_i$ = ($\pi = (x_1\;x_2))\gamma_i$
\State Initialize $down = 0$, $up = 0$
\ForAll {$\alpha \in $ set of fixed points in  $\pi_i D_n $}
  \If{$(\alpha \preceq (\beta \ | \ \gamma))$} \Comment{"|" is bitwise "OR"}
    \State $down$ += $1$
    \EndIf
\EndFor
\ForAll {$\delta \in $ set of fixed points in  $\pi_i D_n $}
  \If{$((\beta \ \& \ \gamma) \preceq \delta)$} \Comment{"\&" is bitwise "AND"}
    \State $up$ += $1$
    \EndIf
\EndFor
\State $k$ += $up \times down$
\EndFor
\end{algorithmic}
\end{algorithm}

\section{Implementation and results}

The algorithms have been implemented in Java and run on a computer with an Intel Core i7-9750H processor. Algorithms have been tested and compared with the results of Chuchang and Shoben \cite{Ch1986} for $r_7$ where some misprints have been found, therefore I give a complete, correct table of detailed calculation results for $r_7$. The total computation time of $r_8$ was approximately a few minutes (with $d_8$ precomputed).

\begin{table}[H]
\centering
    {\renewcommand{\arraystretch}{1.1}
    \begin{tabular}{c*{3}{c}c}
    $i$              & $\pi_i$ & $\mu_i$ & $\phi(\pi_i)$ \\
    \hline
        1    & (1)              & 1     & 2414682040998  \\
        2    & (12)             & 21    & 2208001624  \\
        3    & (123)            & 70   & 2068224  \\
        4    & (1234)           & 210   & 60312  \\
        5    & (12345)          & 504  & 1548  \\
        6    & (123456)         & 840  & 766  \\
        7    & (1234567)        & 720  & 101  \\
        8    & (12)(34)         & 105   & 67922470  \\
        9    & (12)(345)        & 420  & 59542  \\
        10   & (12)(3456)       & 630  & 26878  \\
        11   & (12)(34567)      & 504  & 264  \\
        12   & (123)(456)       & 280  & 69264  \\
        13   & (123)(4567)      & 420  & 294  \\
        14   & (12)(34)(56)     & 105   & 12015832  \\
        15   & (12)(34)(567)    & 210  & 10192  \\
        
        \hline
    \end{tabular}}
    \[ r_7 = \frac{1}{5040} \sum_{i=1}^{k = 15} \mu_i \phi(\pi_i) = 490013148 \]
    \vspace{-2mm}
    \caption{Detailed calculation results for $r_7$} \label{tab:r7results}
    \end{table}

\vspace{6mm}

\begin{table}[H]
\centering
    {\renewcommand{\arraystretch}{1.1}
    \begin{tabular}{c*{3}{c}c}
    $i$              & $\pi_i$ & $\mu_i$ & $\phi(\pi_i)$ \\
    \hline
        1    & (1)              & 1     & 56130437228687557907788  \\
        2    & (12)             & 28    & 101627867809333596  \\
        3    & (123)            & 112   & 262808891710  \\
        4    & (1234)           & 420   & 424234996  \\
        5    & (12345)          & 1344  & 531708  \\
        6    & (123456)         & 3360  & 144320  \\
        7    & (1234567)        & 5760  & 3858  \\
        8    & (12345678)       & 5040  & 2364  \\
        9    & (12)(34)         & 210   & 182755441509724  \\
        10   & (12)(345)        & 1120  & 401622018  \\
        11   & (12)(3456)       & 2520  & 93994196  \\
        12   & (12)(34567)      & 4032  & 21216  \\
        13   & (12)(345678)     & 3360  & 70096  \\
        14   & (123)(456)       & 1120  & 535426780  \\
        15   & (123)(4567)      & 3360  & 25168  \\
        16   & (123)(45678)     & 2688  & 870  \\
        17   & (1234)(5678)     & 1260  & 3211276  \\
        18   & (12)(34)(56)     & 420   & 7377670895900  \\
        19   & (12)(34)(567)    & 1680  & 16380370  \\
        20   & (12)(34)(5678)   & 1260  & 37834164  \\
        21   & (12)(345)(678)   & 1120  & 3607596  \\
        22   & (12)(34)(56)(78) & 105   & 2038188253420  \\
        
        \hline
    \end{tabular}}
    \[ r_8 = \frac{1}{40320} \sum_{i=1}^{k = 22} \mu_i \phi(\pi_i) = 1392195548889993358 \]
    \vspace{-2mm}
    \caption{Detailed calculation results for $r_8$} \label{tab:r8results}
    \end{table}
    
    \newpage

\bibliographystyle{unsrt}

\begin{thebibliography}{1}

\bibitem{Ch1985}
Liu Chuchang, Hu Shoben.
\newblock A mechanical algorithm of equivalent classification for free distributive lattices.
\newblock in {\em Chinese Journal of Computers, 1985, Issue 02} (in Chinese)

\bibitem{Ch1986}
Liu Chuchang, Hu Shoben.
\newblock A note on computation of the numbers of equivalent classes for the free distributive lattices
\newblock in {\em Journal of Wuhan University (Natural Science Edition), 1986, Issue 01}, pages 13-17 (in Chinese)

\bibitem{r7}
Tamon Stephen, Timothy Yusun.
\newblock Counting inequivalent monotone Boolean functions
\newblock in {\em Discrete Applied Mathematics, Volume 167, 20 April 2014}, pages 15-24

\bibitem{Wied1991}
Doug Wiedemann.
\newblock A computation of the eighth Dedekind number
\newblock in {\em Order 8, 1991}, pages 5-6

\bibitem{Fid2001}
Robert Fidytek et al. 
\newblock Algorithms counting monotone Boolean functions
\newblock in {\em Information Processing Letters, 2001, issue 79}, pages 203-209

\bibitem{Bakoev}
Valentin Bakoev. 
\newblock One more way for counting monotone Boolean
functions
\newblock in {\em Thirteenth International Workshop on Algebraic and Combinatorial Coding Theory, 2012}, pages  47–52

\bibitem{earl2014}
Richard Earl.
\newblock Groups and group actions
\newblock {\url{https://courses.maths.ox.ac.uk/node/view_material/43836}}, accessed 01.05.2021

\bibitem{assarpour}
Ali Assarpour.
\newblock List, Sample, and Count (2018). CUNY Academic Works.
\newblock {\url{https://academicworks.cuny.edu/gc_etds/2869}}, accessed 01.05.2021

\bibitem{campo}
Frank a Campo.
\newblock Relations between Powers of Dedekind Numbers and Exponential Sums Related to Them
\newblock in {\em Journal of Integer Sequences, Vol. 21 (2018)}

\bibitem{oeis}
Neil J.A. Sloane.
\newblock The Online Encyclopedia of Integer Sequences
\newblock {\url{https://oeis.org}}

\end{thebibliography}

\end{document}